\documentclass{amsart}

\usepackage{amsmath,amssymb}

\usepackage[dvips]{graphicx}
\usepackage{pstricks}

\newcommand{\dR}{{\mathbb{R}}}
\newcommand{\sL}{{\mathcal{L}}}
\newcommand{\beq}[2]{\begin{equation}\label{#1}#2} %follow with \end{equation}}

\DeclareMathOperator{\cov}{cov}

\newtheorem{thm}{Theorem}
\newtheorem{defn}[thm]{Definition}
\newtheorem{exmp}[thm]{Example}
\newtheorem{lem}[thm]{Lemma}
\newtheorem{rem}[thm]{Remark}
\newtheorem{prop}[thm]{Proposition}

\begin{document}

\title{Densities for Random Balanced Sampling}

\author{Peter Bubenik}
\address{Department of Mathematics, Cleveland State University, 
    2121 Euclid Ave. RT 1515, Cleveland OH 44115, USA}
\email{p.bubenik@csuohio.edu}
\urladdr{http://academic.csuohio.edu/bubenik\_p/}
\thanks{The first author
    thanks NSERC of Canada for support under its Summer Undergraduate
    Scholarship program; his work during summer 1995 was the basis for
    some of the results reported here.}

\author{John Holbrook}
\address{Department of Mathematics and Statistics,
University of Guelph, Guelph, Ontario, Canada N1G2W1}
\email{jholbroo@uoguelph.ca}
\thanks{The second author also thanks NSERC
    for support through its Discovery Grants program.}

\begin{abstract}
A random balanced sample (RBS) is a multivariate distribution with $n$ components $X_1,\ldots,X_n$, each uniformly distributed 
on $[-1,1]$, such that the sum of these components is precisely 0. The corresponding vectors $\vec X$ lie in an
$(n-1)$--dimensional polytope $M(n)$. We present new methods for the
construction of  such RBS via densities over $M(n)$ and these apply
for arbitrary $n$. While simple densities had been known previously
for small values of $n$ (namely 2, 3 and 4), for larger $n$ the known distributions
with large support were fractal distributions (with fractal dimension
asymptotic to $n$ as $n\to\infty$). Applications of RBS distributions
include sampling with antithetic coupling to reduce variance, and the
isolation of nonlinearities. We also show
that the previously known densities (for $n\leq4$) are in fact the
only solutions in a natural and very large class of potential RBS
densities. This finding clarifies the need for new methods, such as those presented here.
\end{abstract}


%AMS Subject Classification: 

\maketitle

\section{ Introduction} \label{section:1}

The aim of this article is to improve our understanding of ``random balanced samples'' (RBS), multivariate distributions having specified marginals as well as a ``balanced'' property; the precise definition is given below. Random balanced samples have found several applications, described in the literature, but the construction of nontrivial RBS distributions has presented interesting mathematical challenges. We describe new methods that generate RBS \emph{densities} with respect to the underlying Lebesgue measure. In contrast, previous work on this problem has relied on fractal geometry and iterated function systems. Ironically, the earliest examples of RBS distributions were given by simple explicit densities. In Section~\ref{section:3}, however, we show that those earlier techniques \emph{cannot} be extended to higher dimensions.
\begin{defn}
A \emph{random balanced sample} (RBS) of size $n$ is a system of random variables 
\[
    \vec X = (X_1,X_2,\dots,X_n)
\]
such that $\vec X$ is \emph{balanced}, ie
\[
     \sum_{k=1}^n X_k = 0,
\]
and each $X_k$ is uniformly distributed over $[-1,1]$. Other marginal
distributions might be considered, but we focus here on the uniform
case.
\end{defn}

Our main concern here is to present new methods for constructing RBS densities (see section 2), and to point out the severe limitations of older methods (see section 3). Nevertheless, to set the context we briefly review some situations where RBS distributions may be usefully applied.

A classical technique for reducing Monte Carlo variance is the use of ``antithetic variates''.
First introduced by Hammersley and Morton in 1956~\cite{hammersleyMorton}, it has been extended to larger groupings of variables by Arvidsen and Johnsson~\cite{arvidsenJohnsson}, and applied to the bootstrap by P. Hall~\cite{hall:antithetic}.
For more recent applications, see the papers of Craiu and Meng \cite{craiuMeng:multiprocessPAC,craiuMeng:antitheticCouplingFPS,craiuMeng:chanceAndFractals}.
In the following example we motivate the use of antithetic variates, and we show that in the exchangeable case, RBS distributions are \emph{extremely antithetic}. 
That is, the variables are as negatively correlated as possible.

\begin{exmp} \label{exmp:1.1}
We may wish to estimate the mean $\frac{1}{2}\int_{-1}^1 f(x)\,dx$ of an unknown function $f:[-1,1]\to\dR$ by means of 
\[
   \overline{f}_n=\frac1n\sum_{k=1}^n f(X_k).
\]
The variance of $\overline{f}_n$ will be reduced if $\cov(f(X_k),f(X_j))<0$ for $k\neq j$. 
This reduction can be achieved for a variety of functions $f$ by means of \emph{antithetic coupling} of the sampling variables $X_k$, ie by insisting that $\cov(X_k,X_j)$ be negative for $k\neq j$. Note that for $X_k$ uniformly distributed over $[-1,1]$ we have
\[
    E\left(\left|\sum_{k=1}^n X_k\right|^2\right)=\frac{n}3+\sum_{k\neq j}E(X_kX_j),
\]
so that if $\cov(X_k,X_j)=\alpha$ for all $k\neq j$ (eg in the exchangeable case) then
\[
    \alpha=\frac1{n(n-1)}\left(E\left(\left|\sum_{k=1}^n X_k\right|^2\right)-\frac{n}3\right).
\]
Thus the minimum value of $\alpha$ is $-\frac{1}{3(n-1)}$, achieved precisely
when we have a RBS.
\end{exmp}

\begin{exmp} \label{exmp:1.2}
In the same setting as Example~\ref{exmp:1.1}, the estimate $\overline{f}_n$ via 
a RBS allows a cleaner distinction between the linear and nonlinear
parts of $f$. Indeed, if $f(x)=Cx+g(x)$ where $g$ is nonlinear then a
RBS entirely eliminates the linear component $Cx$ (even though the
constant $C$ may not be known) from the estimate $\overline{f}_n$. Applications of this nature are discussed in Gerow et al. ~\cite{gerowMcCulloch,gerowHolbrook:mathIntelligencer,gerow:phdThesis,gerow:mastersThesis}; see also ~\cite{royallCumberland}.
\end{exmp}

The structure of RBS distributions is by no means a simple matter; there is a bewildering variety of RBS distributions for sample size $n\geq3$. Each may be regarded as a probability distribution on the regular polytope
\[
    M(n) = \left\{\vec x \in [-1,1]^n:\sum_{k=1}^n x_k=0\right\},
\]
with the additional requirement of marginal uniformity. For small values of $n$, the polytopes $M(n)$ are familiar geometric objects; for example, one easily sees that $M(3)$ is the regular hexagon with vertices $(\pm1,\mp1,0),(\pm1,0,\mp1)$, and $(0,\pm1,\mp1)$. It seems that
the earliest solution to this puzzle (aside from the trivial case
where $n=2$) is due to D. Robson. He noted that a simple piecewise
linear density on $M(3)$ meets the requirements. See
Figure~\ref{figure:1} below. K. Gerow, in~\cite{gerow:mastersThesis}, extended Robson's method to $n=4$, and detected problems for this method in case $n>4$. As it turned out, even a broad generalization of the Gerow--Robson method could not produce RBS densities for $n\geq5$. This phenomenon is studied in detail in \S\ref{section:3} below.

\begin{figure} 
\begin{center}
\includegraphics[width=8cm,keepaspectratio=true]{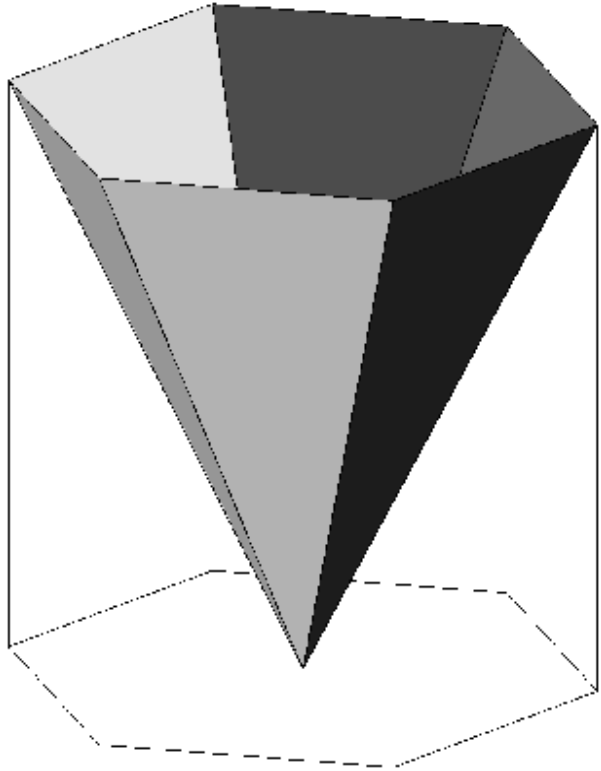}
\end{center}
\caption{The Robson density graphed as a piecewise linear density over
  $M(3)$, represented as a regular hexagon. Gerow discovered that a
  similar construction (density proportional to $\max|X_k|$) works on
  $M(4)$.}
\label{figure:1}
\end{figure}

As alternatives to the Robson density on $M(3)$, many other RBS
distributions were constructed; part of~\cite{gerowHolbrook:mathIntelligencer} is devoted
to a survey of these, many of which have attractive
geometric features. In particular, certain \emph{fractal} measures
define RBS distributions for arbitrary sample size $n$. See Figures
\ref{figure:2} and \ref{figure:3} below. While such fractal measures
cannot have all of $M(n)$ as their support, the fractal dimension of
their supports can be made large with respect to the topological
dimension $n-1$ of $M(n)$ (see \cite{gerowHolbrook:mathIntelligencer,ramlochan:thesis}). For some time, there appeared to be a dichotomy in the possible RBS constructions: for small $n$ we had simple RBS densities, but for larger $n$ only \emph{singular} RBS  distributions (the fractal constructions, which produce measures singular with respect to the natural ($n-1$)--volume on $M(n)$). 

\begin{figure} 
\begin{center}
\includegraphics[width=8cm,keepaspectratio=true]{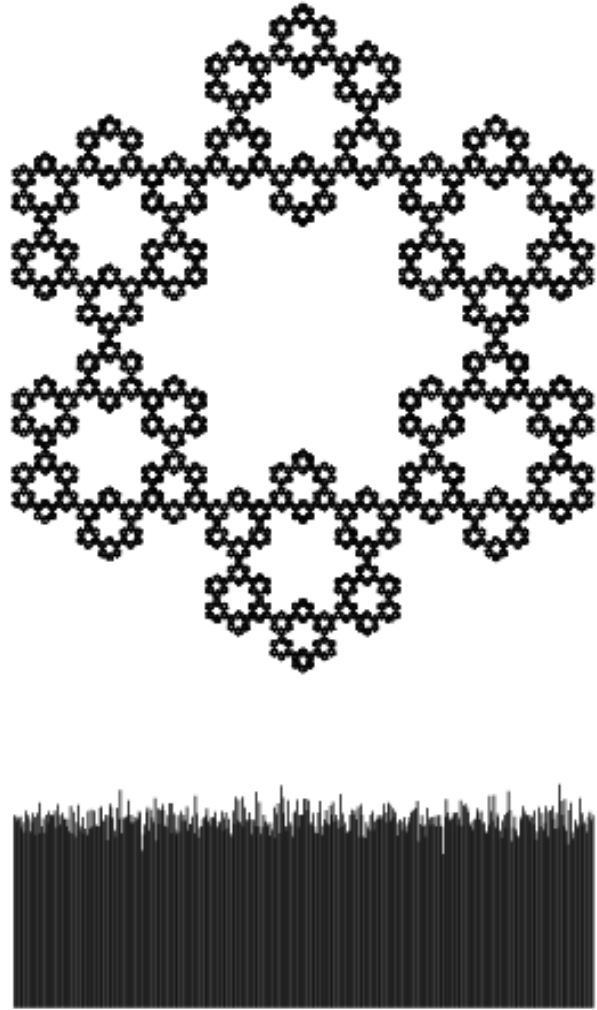}
\end{center}
\caption{A fractal ``superstar'' RBS distribution on the hexagon $M(3)$. The histogram generated by this simulation is shown at the bottom, visually verifying the uniform distribution of $X_k$.}
\label{figure:2}
\end{figure}

\begin{figure}
\begin{center}
\includegraphics[width=8cm,keepaspectratio=true]{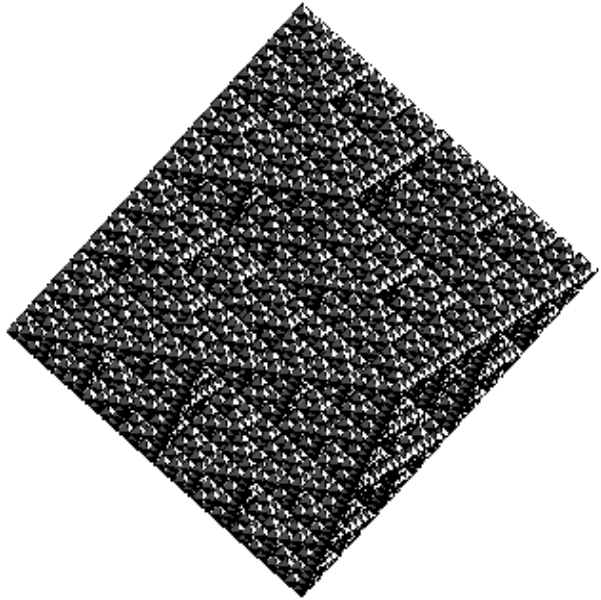}
\end{center}
\caption{As in Figure~\ref{figure:2}, a fractal RBS may be generated in the octahedron $M(4)$ by means of an IFS (iterated function system) mapping all of $M(n)$ closer to its various vertices (6 in this case).}
\label{figure:3}
\end{figure}

More recently, however, we have found constructions that yield RBS densities for any sample size $n$. These are described in \S\ref{section:2} below. They introduce a technique of \emph{redistribution} that may have other applications as well. Like the earlier fractal constructions (which can be implemented via iterated function systems), these new methods are algorithmically effective. That is, workers in the field can easily generate sampling vectors $\vec X$ by following the procedures outlined in \S\ref{section:2}.

\section{Optimal densities} \label{section:2}

\subsection{Goals} \label{subsection:2.1}

Given a RBS, it may be optimized in several directions, and these goals are
to some degree incompatible. We may wish, for example, to achieve as much mutual
independence as possible among the coordinate random variables $X_k$. In this way we would stay near the more familiar situation of i.i.d. sampling. It is clear, however, that for a
given value of $n$ not more than $n/2$ of the $X_k$ can be mutually
independent. Suppose, for example, that $X_1,X_2,\dots,X_m$ are mutually independent;
then there is a positive probability (namely $(\epsilon/2)^m$, if $\epsilon>0$ is small)
that
\begin{equation}\label{E:2.1}
     X_k\in[1-\epsilon,1]\quad(k=1,\dots,m).
\end{equation}
In this case we would have
\begin{equation}\label{E:2.2}
   m(1-\epsilon)\leq\sum_{k=1}^m X_k = -\sum_{k=m+1}^n X_k\leq (n-m),
\end{equation}
which is not possible (for sufficiently small $\epsilon$) if $2m>n$.
On the other hand, if we compromise in other directions, the maximal mutual
independence is easy to obtain. Consider, for example, an even sample size $n=2m$
and a RBS defined by choosing $X_1,X_2,\dots,X_m$ independently (and each $X_k$ uniform
in $[-1,1]$, of course), then setting
\begin{equation}\label{E:2.3}
       X_{m+k} = -X_k \quad (k=1,\dots,m).
\end{equation}
Here we have the coordinates partitioned into two subsets of $m=n/2$
mutually independent random variables (for example, the first $m$ and the last $m$
of the $X_k$). This RBS, however, is \emph{degenerate} in the sense that the distribution
of $\vec{X}$ is supported on a $(n/2)$--dimensional subset of the natural range $M(n)$
for RBS of size $n$. Recall that
\begin{equation}\label{E:2.4}
   M(n) = \left\{\vec{x}\in [-1,1]^n:\sum_{k=1}^n x_k = 0\right\},
\end{equation}
so that $M(n)$ has the much larger dimension $n-1$. A sampling procedure based on such a RBS lacks robustness in a certain sense: it may overlook significant structure in the distribution under investigation because that structure lies outside the support of the RBS. Following this line of thought, the  dimension of the support of a RBS
distribution has been viewed as a measure of \emph{robustness}; in some interesting
cases this dimension must be computed as a \emph{fractal} dimension
(see~\cite{gerowHolbrook:mathIntelligencer}, for example; a fractal construction is also displayed in \cite{craiuMeng:chanceAndFractals}). In \S\ref{subsection:2.2} we shall see how the degenerate RBS defined by (\ref{E:2.3})
can be modified to achieve maximal robustness $(n-1)$ without much loss in mutual
independence of the coordinates.

A natural $(n-1)$--dimensional model for $M(n)$ is obtained by first
constructing $n$ unit vectors $u_1,u_2,\dots,u_n$ in $\dR^{n-1}$ such
that $(u_i,u_j)=-1/(n-1)$ whenever $i\neq j$. Then $M(n)$ may be
identified with $\{v\in \dR^{n-1}:(\forall
k)\,\,-1\leq(v,u_k)\leq1\}$. Corresponding to each such $v$ we have
$\vec X$ with $X_k=(v,u_k)$ ($k=1,2,\dots,n$). See~\cite{gerowHolbrook:mathIntelligencer} for
further details. These models of $M(n)$ yield regular polytopes in
$\dR^{n-1}$. For example, $M(3)$ is seen as a regular hexagon, as in
Figures \ref{figure:1} and \ref{figure:2}, while $M(4)$ is the regular
octahedron outlined in Figure~\ref{figure:3}. Figure~\ref{figure:4}
below presents one of the 3--dimensional faces of $M(5)$, which is
itself a 4--dimensional object. 

\begin{figure} 
\begin{center}
\includegraphics[width=8cm,keepaspectratio=true]{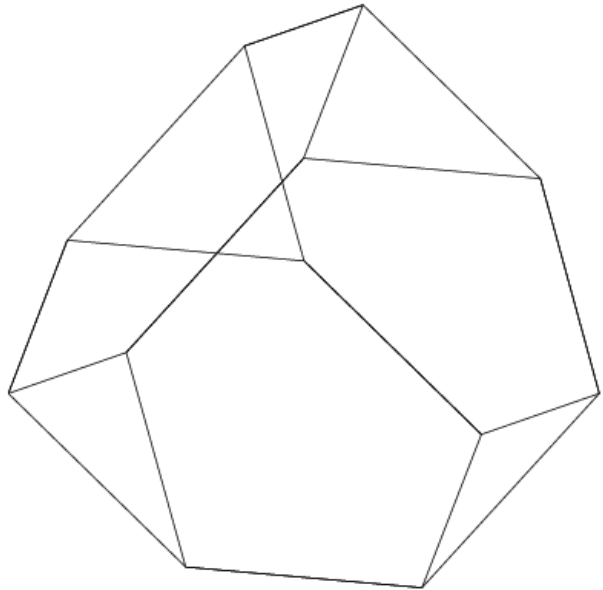}
\end{center}
\caption{While the 4--dimensional $M(5)$ cannot be easily visualized, each of its faces is similar to the elegant 3--dimensional object seen here. It may be regarded as a truncated tetrahedron; the result has 4 regular hexagons and 4 equilateral triangles as \emph{its} faces.}
\label{figure:4}
\end{figure}

Along with independence and robustness, a third desirable feature of a RBS procedure
might be \emph{algorithmic efficiency}, the efficiency with which sampling variables 
may be generated for use in experiments. The RBS defined by \eqref{E:2.3}, for example,
is very efficient, since it involves very little beyond $m$ calls to a random number
generator to produce a sampling vector $\vec{X}$. We shall see that the modifications
proposed in \S\ref{subsection:2.2} retain most of this efficiency as they increase robustness.

\subsection{Redistribution} \label{subsection:2.2}

Here we illustrate, in a simple setting, a general procedure of 
\emph{redistribution} that may be used to increase the robustness of a RBS. 
In this initial example, we confine ourselves to the redistribution of \emph{pairs}
of independent coordinates chosen from $\vec{X}$.

\begin{lem} \label{lem:2.1}
Given independent $X_1,X_2$, each 
uniformly distributed on $[-1,1]$,
let $S=X_1+X_2$ and define the new variables 
\begin{equation}\label{E:2.5}
     Y_1=\frac{S}2+\left(1-\frac{|S|}2\right) T,\quad
     Y_2=\frac{S}2-\left(1-\frac{|S|}2\right) T,
\end{equation}
where $T$ is uniform on $[-1,1]$ and independent of the $X_k$; then $Y_1,Y_2$ are also
independent and uniform on $[-1,1]$.
\end{lem}

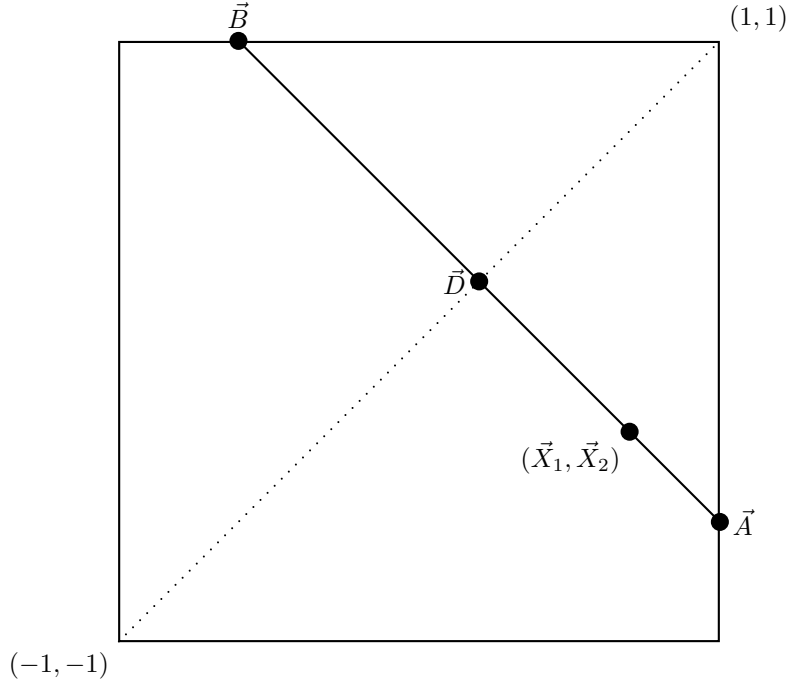
\begin{figure}
\psset{unit=4cm}
\centering
\begin{pspicture}(-1.5,-1.5)(1.5,1.5)
\psframe(-1,-1)(1,1)
\psline[linestyle=dotted](-1,-1)(1,1)
\psline(1,-0.6)(-0.6,1)
\qdisk(1,-0.6){0.03}
\qdisk(-0.6,1){0.03}
\qdisk(0.2,0.2){0.03}
\qdisk(0.7,-0.3){0.03}
\uput[r](1,-0.6){$\vec A$}
\uput[u](-0.6,1){$\vec B$}
\uput[l](0.2,0.2){$\vec D$}
\uput[dl](0.7,-0.3){$(\vec X_1,\vec X_2)$}
\uput[dl](-1,-1){$(-1,-1)$}
\uput[ur](1,1){$(1,1)$}
\end{pspicture}%
\caption{A graphical view of the redistribution algorithm for the case $\vec X_1 = 0.7$, $\vec X_2=-0.3$.} \label{figRedistribution} 
\end{figure}

\begin{rem} 
Since $Y_1+Y_2=S=X_1+X_2$, we may think of 
the procedure in this lemma
as \emph{redistributing} the part of a RBS captured by $X_1+X_2$. The geometry behind this construction is revealed by writing $\vec Y=(Y_1,Y_2)=\vec D + T\vec Q$ where $\vec D=(S/2,S/2)$ is a point on the diagonal of $[-1,1]^2$ and 
\[
     \vec Q=\left(1-\frac{|S|}{2},-\left(1-\frac{|S|}{2}\right)\right), 
\]
so that $\vec A=\vec D+\vec Q$ and $\vec B=\vec D-\vec Q$ are the endpoints of the segment in $[-1,1]^2$ passing through $\vec D$ and perpendicular to the diagonal
(see Figure~\ref{figRedistribution}).
Thus $\vec Y$ is chosen in two steps: first $\vec D$ is chosen with a density proportional to the length of $\vec A\vec B$, then $\vec Y$ is placed at a point chosen uniformly along $\vec A\vec B$. This procedure strongly suggests that $\vec Y$ will be uniformly distributed over the square $[-1,1]^2$. We present a more formal proof below, computing the joint density of $(Y_1,Y_2)$.
\end{rem}

\begin{proof}
Let $D=S/2$. The joint density $f(t,d)$ of $(T,D)$ on $[-1,1]^2$ is given by $\frac{1}{2}(1-|d|)$, since $T$ and $D$ are independent and it is easy to compute the density of $D$ on $[-1,1]$ as $1-|d|$. Also
\[
    (Y_1,Y_2)=(D+(1-|D|)T,D-(1-|D|)T),
\]
a transformation with Jacobian matrix
\[
J(t,d)=
\begin{bmatrix} (1-|d|)&1-\mbox{sign}(d)t\\
                 -(1-|d|)&1+\mbox{sign}(d)t
\end{bmatrix}
\]
(compute for the two cases $d\geq0$ and $d<0$). It follows that $\det J(t,d)=2(1-|d|)$. The transformation is injective from $[-1,1]^2$ to $[-1,1]^2$, so that we may compute (compare the discussion of ``change of variable'' in ~\cite{Hoel54}, page 213) the joint density of $(Y_1,Y_2)$ at the image of $(T,D)$ as $f(t,d)/ \det J(t,d)=1/4$, the uniform density on $[-1,1]^2$. Thus $Y_1,Y_2$ are uniform over $[-1,1]$ and are independent.
\qed
\end{proof}

The following theorem shows how to redistribute selected pairs of
coordinates from the degenerate RBS defined by \eqref{E:2.3} so that
we obtain a RBS that is maximally 
robust.

\begin{thm} \label{thm:2.2}
Let $X_1,\dots,X_m,T_1,\dots,T_m$ be mutually independent random
variables, each uniform on $[-1,1]$. Assume $m>1$ and let 
\begin{equation}\label{E:2.6}
    S_k = X_k - X_{k+1} \quad (k=1,\dots,m),
\end{equation}
with the understanding that $X_{m+1}=X_1$. Next (redistributing the coordinates combined
in $S_k$ as in the lemma) let
\begin{equation}\label{E:2.7}
   Y_{2k-1}=\frac{S_k}2+\left(1-\frac{|S_k|}2\right) T_k,\quad
   Y_{2k}=\frac{S_k}2 - \left(1-\frac{|S_k|}2\right) T_k,
\end{equation}
for $k=1,\dots,m$. Then $\vec{Y}$ defines a RBS of size $n=2m$ that is
maximally robust, ie that has a set of dimension $n-1$ as its support. Moreover, the distribution
of $\vec Y$ is given by a density on $M(n)$.
\end{thm}

\begin{proof} 
Since $X_k$ and $-X_{k+1}$ are independent and uniform
on $[-1,1]$, Lemma~\ref{lem:2.1} (with appropriate change of notation) ensures
  that each coordinate\footnote{We remark that it may be of interest to assess 
the degree to which these coordinates are mutually independent.} 
of $\vec{Y}$ is uniform on $[-1,1]$.
Clearly $\sum_1^n Y_k = \sum_1^m S_k = 0$, so that
$\vec{Y}$ indeed defines a RBS. To see that it is maximally robust, note first that every
$\vec{y}\in M(n)$ with sufficiently small $|y_k|$ can occur as a value of $\vec{Y}$. Indeed,
if $|y_k|<\frac{1}{n}$ for each $k$, we set (for $k=1,2,\dots,m$) 
\begin{equation}\label{E:2.8}
    x_k=-\sum_{j=1}^{k-1} r_j,
\end{equation}
where $r_j=y_{2j-1}+y_{2j}$.
Note that $x_1=x_{m+1}=0$ and that $|x_k|\leq \frac{2(m-1)}{n}<1$. We have $x_k-x_{k+1}=r_k$ for each $k$; for $k=m$ this
is a consequence of the fact that $\sum_1^m r_k=\sum_1^n y_k=0$.
Solving for appropriate values of $t_k$, we find that we require only that
\[
    t_k=\frac{y_{2k-1}-y_{2k}}{2-|y_{2k-1}+y_{2k}|}.
\]
It is easy to check that $|t_k|\leq1$, simply because $|y_{2k-1}|,|y_{2k}|\leq1$; if $y_{2k-1}=y_{2k}=\pm 1$, any value of $t_k$ will do. Thus, with $\vec X=\vec x$ and $\vec T=\vec t$, we have $\vec Y=\vec y$, as claimed.

Next consider a neighborhood $\vec y(\epsilon)$
of such a $\vec{y}$ in $M(n)$; we must show that the procedure of the theorem places sampling
vectors in $\vec y(\epsilon)$ with positive probability. By the continuity of the procedure for obtaining $\vec Y$ from $(\vec X,\vec T)$,
there is a neighborhood $\vec u(\delta)$ of $(\vec{x},\vec{t})$ in $[-1,1]^n$ that is mapped
by the procedure into $\vec y(\epsilon)$. The probability of $\vec u(\delta)$ is already positive;
indeed it's just the
normalized $n$--volume of $\vec u(\delta)$, since the coordinates  of $(\vec X,\vec T)$ are chosen independently and uniformly 
in $[-1,1]$.

To see that the distribution of $\vec Y$ is given by a density on $M(n)$, we introduce the (linear) mapping $L:[-1,1]^m\to 2M(m)$ defined by
\begin{equation}\label{E:2.9}
    L(x_1,\dots,x_m) = (x_1-x_2,x_2-x_3,\dots,x_m-x_1).
\end{equation}
The analysis above shows that $\vec y\in M(n)$ occurs as a value of $\vec Y$ exactly when $L(\vec x)=\vec r$ for some $\vec x\in [-1,1]^m$. Given $\epsilon>0$, consider the map 
\[
     f:(1-\epsilon)M(n)\to 2M(m)\times[-1,1]^m
\]
defined by $f(\vec y)=(f_1(\vec y),f_2(\vec y))$, where $f_1(\vec y)=\vec r$ and $f_2(\vec y)=\vec t$.
The denominators in the expressions for the $t_k$ are bounded away from zero since $2-|y_{2k-1}+y_{2k}|\geq 2 - 2(1-\epsilon) = 2\epsilon$, so that $f$ is Lipschitz on $(1-\epsilon)M(n)$ for each fixed $\epsilon$. Given a Borel subset, $B$, of $(1-\epsilon)M(n)$, the probability assigned to $B$ by the distribution of $\vec Y$ is the (normalized) $n$--volume (or Lebesgue measure) of $\{(\vec x,\vec t):(L(\vec x),\vec t)\in f(B)\}$. Since $L$ is linear and $[-1,1]^n$ is bounded, this is at most a constant times the $(n-1)$--volume of $f(B)$. As $f$ is Lipschitz, this is in turn bounded by a constant times the $(n-1)$--volume of $B$. Certainly, then, the distribution is absolutely continuous with respect to normalized Lebesgue measure on $(1-\epsilon)M(n)$. Considering a sequence of $\epsilon$--values tending to 0, we see that the distribution is absolutely continuous on $M(n)$ itself and so given by a density with respect to the normalized $(n-1)$--volume on $M(n)$.
\qed
\end{proof}

\subsection{Using all of $M(n)$} \label{subsection:2.3}

The procedure of the last section retains algorithmic efficiency and it is
robust in the sense of dimension, but in most cases the sampling values do not fill all of
$M(n)$. Here we show how to modify the construction to obtain  RBS 
procedures that have \emph{all} of $M(n)$ as support. First let us clarify the
reasons for the failure of the construction in \S\ref{subsection:2.2} to cover $M(n)$; we reuse 
the notations of the proof above. As $\vec{y}$ ranges over $M(n)$
the corresponding $\vec{r}$ fills all of $2M(m)$ (recall that $n=2m$). The procedure
will yield every $\vec{y}\in M(n)$ as a possible value exactly when 
$2M(m)=L([-1,1]^m)$.
Since $L([-1,1]^m)$ is the convex hull of the images of the $2^m$ extreme points
of $[-1,1]^m$, we can simply check whether these images include the extreme points
of $2M(m)$, which are easy to identify (at most one coordinate can
differ from $\pm2$).

For $m=2$, for example, the extreme points of $2M(2)$ are
\begin{equation}\label{E:2.10}
   (2,-2)=L(1,-1) \mbox{ and } (-2,2)=L(-1,1),
\end{equation}
so that the procedure of \S\ref{subsection:2.2} for $m=2$ gives a RBS distribution on $M(4)$ having
a density that is positive everywhere on the octahedron. One can check that this density
is unbounded, in contrast to the piecewise linear density found by Gerow (see section 3,
below). For $m=3$ we again get a RBS distribution with all of $M(6)$ as support, because the
extreme points of $2M(3)$ are:
\begin{equation}\label{E:2.11}
  \pm(2,-2,0)=L(\pm(1,-1,1)),\pm(0,2,-2)=L(\pm(1,1,-1)),\mbox{ etc.}
\end{equation}

For $m=4$, however, one can compute in a similar fashion that $L([-1,1]^4)$ is 
a sort of rhombic dodecahedron lying \emph{strictly} inside the octahedron 
$2M(4)$. See Figure~\ref{figure:5} below. Thus, except in a few simple cases, the procedure of \S\ref{subsection:2.2}, while it is
maximally robust in terms of dimension, does not fill $M(n)$. As we'll now see,
the remedy is simply to symmetrize this procedure.

\begin{figure} 
\begin{center}
\includegraphics[width=8cm,keepaspectratio=true]{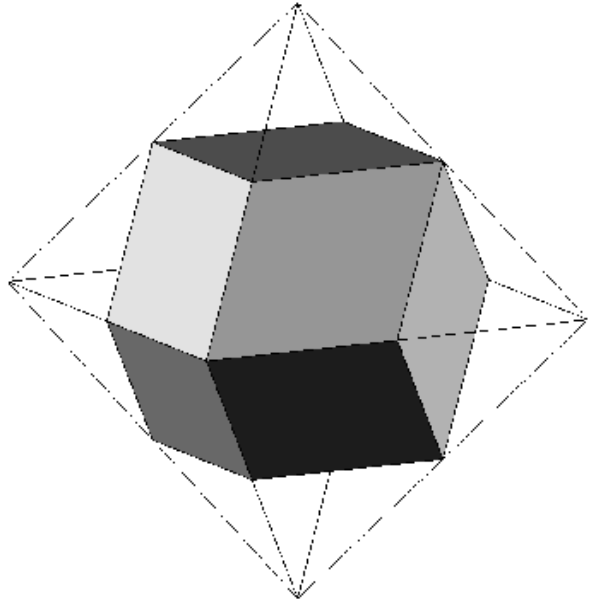}
\end{center}
\caption{The fact that the rhombic dodecahedron does not fill up all of the octahedron $2M(4)$ reveals the need for symmetrization in order to generate a RBS supported on \emph{all} of $M(8)$.}
\label{figure:5}
\end{figure}

\begin{lem} \label{lem:2.3}
Given any $\vec{w}\in M(n)$, there is a permutation $\sigma$ of
$\{1,2,\dots,n\}$ such that for $\vec{z}=\sigma(\vec{w}):=(w_{\sigma(1)},\dots,w_{\sigma(n)})$
all partial sums
\begin{equation}\label{E:2.12}
    \sum_{j=1}^k z_j \quad (k=1,\dots,n)
\end{equation}
lie in [-1,1].
\end{lem}

\begin{proof}
We can choose the values of $\sigma(i)$ inductively. Suppose that
(distinct) $\sigma(1),\dots,\sigma(i)$ have been chosen so that the partial sums
\eqref{E:2.12} lie in $[-1.1]$ for all $k\leq i$. Since
\begin{equation}\label{E:2.13}
   \sum_{j=1}^i w_{\sigma(j)} + \sum_{j\notin\sigma(\{1,\dots,i\})} w_j = 0,
\end{equation}
we can choose $\sigma(i+1)$ such that $w_{\sigma(i+1)}$ has sign opposite from that
of the first sum in (\ref{E:2.13}); it follows that the partial sum (\ref{E:2.12})
lies in $[-1,1]$ for $k=i+1$ as well.
\qed
\end{proof}

\begin{thm} \label{thm:2.4}
Given $n=2m$, let the RBS $\vec{Y}$ be defined as in 
Theorem~\ref{thm:2.2}, and let $\sigma$ be a randomly chosen permutation of $\{1,\dots,n\}$. Then
$\vec{W}=\sigma(\vec{Y})$ defines a RBS with all of $M(n)$ as support.
\end{thm}

\begin{proof}
For any $\vec{w} \in M(n)$ using the $\sigma$ from Lemma~\ref{lem:2.3},
$\vec{w} = \sigma^{-1}(\vec{z})$ where all the partial sums of 
(\ref{E:2.12}) lie in $[-1,1]$. 
The proof of Theorem~\ref{thm:2.2} showed that any such $\vec{z}$
could occur as a RBS $\vec{Y}$ defined in \S\ref{subsection:2.2}.
Thus $\vec{w}$ can occur as $\tau(\vec{Y})$ for the permutation
$\tau = \sigma^{-1}$.
\qed
\end{proof}

\subsection{Odd sample sizes} \label{subsection:2.4}

For odd sample sizes, $n=2m+1$ we can likewise construct an algorithmically efficient 
degenerate RBS with two disjoint maximal subsets of mutually independent variables. 
This can be done by choosing $X_{1}, \ldots, X_{m}$ (with each $X_{k}$ uniform in $[-1,1]$)
and $B \in \{-1,1\}$ independently, and then setting 
\begin{gather*}
    X_{m+k} =  - X_{k} \quad (k=1, \ldots, m-1) \\
    X_{2m} =  - \frac{1}{2} (X_{m}+B) \\
    X_{2m+1} =  - \frac{1}{2} (X_{m}-B) 
\end{gather*}
It is easy to check that this gives a RBS.
As in the even case, we can make this distribution robust by redistributing pairs of the variables.

\begin{thm} \label{thm:2.5}
    Let $X_{1}, \ldots , X_{m}, T_{1}, \ldots , T_{m}$ be mutually 
    independent random variables, each uniform on $[-1,1]$ and $B$ be a 
    discrete random variable,   uniform on $\{-1,1\}$.
    Assume $m > 1$ and let 
\begin{gather*}
    S_{k} = X_{k} - X_{k+1}, \quad (k=1, \ldots m-2) \\
    S_{m-1} = X_{m-1} - \frac{1}{2} (X_{m} + B) \\
    S_{m} = - \frac{1}{2} (X_{m} - B) - X_{1} 
\end{gather*}
    Next redistribute the variables as in Lemma~\ref{lem:2.1}:
\begin{gather*}
    Y_{2k-1}  = \frac{S_k}2+\left(1-\frac{|S_k|}2\right) T_k, \quad (k=1, \ldots, m) \\
    Y_{2k} = \frac{S_k}2 - \left(1-\frac{|S_k|}2\right) T_k  \quad (k=1, \ldots, m) \\
    Y_{2m+1} = X_{m}
\end{gather*}
    Then $\vec{Y}$ is an RBS of size $n=2m+1$ that is maximally robust.
That is, it has a set of dimension $n-1$ as its support.
\end{thm}

\begin{proof}
    First note that $- \frac{1}{2} (X_{m} + B)$ and $- \frac{1}{2} (X_{m} - B)$ are uniformly distributed
    on $[-1,1]$. As a result, for each $k$, $S_{k}$ is the sum of two independent
    uniform variables, and by Lemma~\ref{lem:2.1}, each coordinate of $\vec{Y}$ is uniform on $[-1,1]$. 
    Also, $\sum_{k=1}^{n} Y_{k} = \sum_{k=1}^{m} S_{k} + Y_{2m+1} = 
    - \frac{1}{2} (X_{m} + B) - \frac{1}{2} (X_{m} - B) + X_{m} = 0$, 
    so $\vec{Y}$ does indeed define an RBS.

    To see that this RBS is maximally robust, we will show that every 
    $\vec{y} \in M(n)$ with sufficiently small coordinates can occur as $\vec{Y}$.
    If $|y_{k}| < \frac{1}{2m}$ for each $k$, construct $\vec{x}$ and $\vec{t}$ as follows.
    Fix $b \in \{-1, 1\}$ and
    let $r_{k} = y_{2k-1} + y_{2k}$  for $k=1,\ldots,m$, and let 
    \[ \nonumber
        x_{m} = - \sum_{i=1}^{m} r_{i} 
    \]
    \[ \nonumber
        x_{k} = \frac{1}{2} \left(x_{m} + b\right) + \sum_{i=k}^{m-1} r_{i}, \quad
            (k=1, \ldots, m-1).
    \]
    Then for $k < m$, $x_{k} = \frac{1}{2} (\sum_{i=1}^{m} (-1)^{a_i} r_i + b)$,
    for some choice of ${a_i}$. 
    So since $|r_k| < \frac{1}{m}$, $|x_k| < 1$ for all $k$. 
    Also for $k \leq m-2$, 
\begin{multline*}  
    s_k = x_k - x_{k+1} = r_{k}, \quad
    s_{m-1} = x_{m-1} - \frac{1}{2} (x_{m} + b) = r_{m-1}, \\ \text{and } 
    s_m = - \frac{1}{2} (x_{m} - b) - x_1 = -x_m - \sum_{i=1}^{m-1} r_i = r_m.
\end{multline*}
    Thus we have a $\vec{x} \in [-1,1]^m$ such that the corresponding $s_k = r_k$
    for each $k$. 
    We can then solve for $t_k \in [-1,1]$ such that 
    $Y_{2k-1} = y_{2k-1}$ and $Y_{2k} = y_{2k}$.
    Also $Y_{2m+1} = x_m = -\sum_{i=1}^{m}r_i = -\sum_{k=1}^{2m}y_k
    = y_{2m+1}$.

    Finally consider the neighborhood $y(\epsilon)$ of $\vec{y}$ in $M(n)$; 
    we must show that the procedure of the theorem places sampling
    vectors in $y(\epsilon)$ with positive probability. By the continuity of the 
    procedure there is a neighborhood $\vec u(\delta)$ of $(\vec{x},\vec{t})$ in 
    $[-1,1]^n$ that is mapped by the procedure into $y(\epsilon)$. 
    The probability of $\vec u(\delta)$ is already positive;
    indeed it's just the normalized $n$--volume, since the coordinates are chosen 
    independently and uniformly in $[-1,1]$.
\qed
\end{proof}

Note that the procedure above will yield a $\vec{y} \in M(n)$ if and only if
$x_k = \sum_{i=2k-1}^{n-3} y_i + \frac{1}{2} (y_n + b)$
lies in $[-1,1]$ for $k=1,\ldots, n-3$. 
By symmetrizing the procedure and proving the following lemma which 
is slightly stronger than necessary, we will get a RBS that fills 
all of $M(n)$.

\begin{lem} \label{lem:2.6}
Given any $\vec{w} \in M(n)$ where n is odd, there exists a permutation
$\sigma$ of $\{1,2, \ldots , n\}$ and a $b \in \{-1, 1\}$
such that for 
$\vec{z} = \sigma(\vec{w}):= (w_{\sigma(1)},\ldots,w_{\sigma(1)})$
all of the sums
\begin{equation} \label{eqnPartialSums} 
   \sum_{i=k}^{n-3} z_i + \frac{1}{2} (z_n + b) \quad (k=1,\ldots, n-3)
\end{equation}
lie in $[-1,1]$.
\end{lem}

\begin{proof}
Choose $\sigma(n)$ such that $z_n = w_{\sigma(n)} \geq 0$. 
Let $b = -1$.
Let $a = -\frac{1}{2} (z_n+b)  \in [0,\frac{1}{2}]$.
The condition now becomes 
\[ S_k := \sum_{i=k}^{n-3} z_i \in [-1+a,1+a] \quad (k=1,\ldots,n-3) \]
Taking $S_{n-2} = 0$ we will prove the above condition by induction.
Assume $\sigma(n-3),\ldots,\sigma(k)$ are defined such that
$S_k \in [-1+a,1+a]$.

If $S_k \in [a, 1+a]$ then 
$\sum_{i=k}^{n-3}w_{\sigma(i)} + w_{\sigma(n)} \geq 0$.
Since 
\[\sum_{i=k}^{n-3}w_{\sigma(i)} + w_{\sigma(n)} +
\sum_{i\notin \sigma(\{k,\ldots,n-3,n\})}w_i = 0 \]
we can choose $\sigma(k-1)$ such that $w_{\sigma(k-1)} \leq 0$.
Then $S_{k-1} \in [-1+a,1+a]$.

If $S_{k-1} \in [-1+a,a]$ then if possible choose $\sigma(k-1)$ such
that $w_{\sigma(k-1)} \geq 0$. 
Then $S_{k-1} \in [-1+a,1+a]$.
Otherwise choose any $\sigma(k-1)$ and 
$\forall i \notin \sigma(\{k-1,\ldots,n-3,n\})$, $w_{\sigma(i)} \leq 0$.
So $\sum_{i=k-1}^{n-3} w_{\sigma(i)} + w_{\sigma(n)} \geq 0$
and thus $S_{k-1} \geq -w_{\sigma(n)} = -1+2a \geq -1+a$.
Therefore $S_{k-1} \in [-1+a,1+a]$.
\qed 
\end{proof}

\begin{thm} \label{thm:2.7}
Given $n=2m+1$, let the RBS $\vec{Y}$ be defined as in the above theorem
and let $\sigma$ be a randomly chosen permutation of $\{1,\ldots,n\}$.
Then $\vec{w} = \sigma(\vec{Y})$ defines a RBS with all of $M(n)$ as 
support.
\end{thm}

\begin{proof}
The proof follows from Theorem~\ref{thm:2.5} and Lemma~\ref{lem:2.6}
in much the same way as Theorem~\ref{thm:2.4} followed from
Theorem~\ref{thm:2.2} and Lemma~\ref{lem:2.3}. 
That is, for any $\vec w \in M(n)$ using the $\sigma$ and $b$ from Lemma~\ref{lem:2.6}, $\vec w = \sigma^{-1}(\vec z)$ where all the sums in \eqref{eqnPartialSums} lie in $[-1,1]$.
The proof of Theorem~\ref{thm:2.5} showed that any such $\vec z$ could occur as a RBS $\vec Y$ as defined in the statement of that theorem. 
Thus $\vec w$ can occur as $\tau(\vec Y)$ for the permutation $\tau = \sigma^{-1}$.
\qed 
\end{proof}

\section{Gerow--Robson densities} \label{section:3}

We are interested in probability densities on $M(n)$ with respect to the measure $v_{n-1}$ defined as follows. 
Let $\vec1$ denote the vector $(1,1,\dots,1)\in\dR^n$, and let $\vec1^\perp$ denote the hyperplane consisting of all vectors perpendicular to $\vec1$. For each Borel subset $S$ of $\vec1^\perp$, let $v_{n-1}(S)$ denote the ($n-1$)--dimensional Lebesgue measure of $S$, regarding $\vec1^\perp$ as an isometric copy of $\dR^{n-1}$. 
Let $V_{n-1}=v_{n-1}(M(n))$.

By a \emph{Gerow--Robson} (G--R) density on $M(n)$ we shall mean a probability density $h$ with respect to $v_{n-1}$ such that 
\beq{3.1}{
     h(\vec X) = f(\|\vec X\|_\infty)\quad(\vec X \in M(n)) }
\end{equation}
for some $f:[0,1]\to[0,\infty)$, where
\beq{3.2}{
       \|\vec X\|_\infty = \max_{1\leq k\leq n}|X_k|. }
\end{equation}
We remark that on $M(3)$ and $M(4)$ the subspaces with constant $\infty$-norm are represented by hexagons and octahedra respectively.
This definition is suggested by the densities discovered by Robson and
Gerow for the cases $n=3,4$. They noted that with $f(s)$ proportional
to $s$, ie with $f_3(s)=C_3s$ and $f_4(s)=C_4s$ for certain constants
$C_k$, the corresponding G--R densities define RBS distributions on
$M(3), M(4)$, respectively (see Figure~\ref{figure:1}). 
These facts will also follow from the
general analysis of G--R densities given below. It was natural to try
to extend this construction to larger values of $n$. We'll prove
below a rather surprising fact: not only do densities with
$f_n(s)=C_ns$ fail to generate RBS distributions when $n>4$, but
(at least for $n\leq 250$)
no other choices of $f_n$ yield RBS distributions. Concisely, 
it seems that
the \emph{only} G--R densities that yield RBS distributions are those
discovered by Robson and Gerow. 

As a historical note, Gerow was already aware via numerical
simulations reported in~\cite{gerow:mastersThesis} that $f(s)=C_5s$ did not yield a RBS
distribution. In~\cite{gerowHolbrook:constructionOfRBS} this was verified theoretically
and it was shown that, in fact, \emph{no} choice of $f_n$ gives a RBS
distribution for $n=5$. 

Here we will extend this result to larger values of $n$.
For $n\geq 6$ we will derive a sufficient condition for the
nonexistence of a G--R density on $M(n)$.
This condition can be verified numerically for $n \leq 250$.
We will also show that the densities given by Gerow and Robson are the
only such densities on $M(3)$ and $M(4)$ and we will reprove the
nonexistence of a G--R density on $M(5)$.

We conjecture that the condition we will derive for the
nonexistence of a G--R density on $M(n)$ holds for all $n \geq 6$.
While we are unable to prove this, our computational results settle
all cases likely to be of practical interest.
For very large samples sizes, constructing random balanced samples is
perhaps less important, since for large values of $n$ a uniform sample from
$[-1,1]^n$ is, in view of the law of large numbers, likely to be very nearly balanced.

For any probability measure $P$ on $M(n)$ we have the corresponding distribution function defined on $[0,1]$ by
\beq{3.3}{
     G(r) = P(rM(n)) = P(\{\vec X\in M(n): \|\vec X\|_\infty\leq r\}). }
\end{equation}
In those cases where $G$ has a density on $[0,1]$ we denote it by $g$:
\beq{3.4}{
     G(r) = P(rM(n)) = \int_0^r g(s)\,ds. }
\end{equation}
If $P$ is defined by a G--R density $h$ corresponding to $f$, we must have
\beq{3.5}{
     G(r) = \int_0^r f(s)(v_{n-1}(sM(n)))'\,ds, }
\end{equation}
and since $v_{n-1}(sM(n))=s^{n-1}v_{n-1}(M(n))=s^{n-1}V_{n-1}$,
\beq{3.6}{
     G(r) = \int_0^r f(s)(n-1)V_{n-1}s^{n-2}\,ds. }
\end{equation}
Thus, $g(s)=(n-1)V_{n-1}s^{n-2}f(s)$ for a G--R density. This relation, along with the fact that $G(1)=1$, allows us to properly normalize $f$. For example, in the cases treated by Robson and Gerow we have $f_n(s)=C_ns$; normalizing $g$ we see that $g_n(s)=ns^{n-1}$. Thus the normalizing constants $C_n$ are given by $\frac{n}{(n-1)V_{n-1}}$. A little computation reveals that $V_2=3\sqrt3$ and $V_3=32/3$, so that, in our models, the densities found by Robson and Gerow correspond to $f_3(s)=s/2\sqrt3$ and $f_4(s)=s/8$.
 
Conversely, given any probability density $g$ on $[0,1]$, we shall see how to construct a probability measure on $M(n)$ with the corresponding G--R density: for $\vec X\in M(n)$,
\beq{3.7}{
     h(\vec X) = f(\|\vec X\|_\infty) = \frac{g(\|\vec X\|_\infty)}
                                             {(n-1)V_{n-1}\|\vec X\|_\infty^{n-2}}. }
\end{equation}
It is natural to ask for which $g$ and for which $n$ we obtain RBS distributions. The answer, as claimed above, is given by the following result.

\begin{thm} \label{thm:3.1}
For $n\leq 250$, the G--R density on $M(n)$ corresponding to a density $g$ on $[0,1]$ defines a RBS distribution in exactly two cases: $n=3$ with $g_3(s)=3s^2$, and $n=4$ with $g_4(s)=4s^3$.
\end{thm}

Before turning to the proof of this theorem, we introduce a specific model for the distribution of $\vec X\in M(n)$ according to the G--R density implied by a density $g$ on $[0,1]$. Consider auxiliary random variables $Y_1$ and $Z_2,\dots,Z_n$ defined as follows: $0\leq Y_1\leq1$ has the given density $g$, and (independently) $(Z_2,\dots,Z_n)$ is uniformly distributed with respect to the ($n-2$)--volume on
\beq{3.8}{
     \{\vec Z\in[-1,1]^{n-1}:\sum_{k=2}^n Z_k = -1\}. }
\end{equation}
Let $Y_k=Y_1Z_k$ for $k=2,\dots,n$, so that $Y_1\geq|Y_k|$ and $\sum_{k=1}^n Y_k=0$. Then $\vec X=\vec Y$ yields sample points in that part of $M(n)$ where $\|\vec X\|_\infty=X_1$. Finally, we symmetrize over $M(n)$:
\beq{3.9}{
     \vec X = \pm (Y_{\sigma(1)},\dots,Y_{\sigma(n)}), }
\end{equation}
where the sign $\pm$ and the permutation $\sigma$ are chosen at random.

\begin{prop} \label{prop:3.2}
Relative to ($n-1$)--volume on $M(n)$, this $\vec X$ has a density of the G--R form (\ref{3.7}).
\end{prop}

\begin{proof}
It suffices to consider a point $\vec X \in M(n)$ where 
\[
   \|\vec X\|_\infty=X_1>|X_k|\quad(k=2,\dots,n),
\]
and the behavior as $\Delta t\to0$ of the probability 
$P(\vec X+\Delta tM(n))$ that a sample point falls in the neighborhood $\vec X+\Delta tM(n)$. Denote this probability, briefly, by $P(\Delta t)$. Then
\[
  P(\Delta t) = P\{\pm=+\}P\{\sigma(1)=1\}
                P\{Y_k\in[X_k-\Delta t,X_k+\Delta t]\quad(k=1,\dots,n)\}. 
\]
For small values of $\Delta t$, $P\{Y_1\in[X_1-\Delta t,X_1+\Delta t]\}$ is close to 
$2\Delta tg(X_1)$, and given $Y_1$ the probability
\[
   P\{Y_k\in[X_k-\Delta t,X_k+\Delta t]\quad(k=2,\dots,n)\}=
\]\[
      P\{Z_k\in[\frac{X_k}{Y_1}-\frac{\Delta t}{Y_1},\frac{X_k}{Y_1}+
      \frac{\Delta t}{Y_1}]\quad(k=2,\dots,n)\}
\]
is nearly proportional to $(\Delta t/X_1)^{n-2}$, since $Y_1\approx X_1$ and the hyperplane $\{\sum_2^nZ_k=-1\}$ makes the same angle with each coordinate axis. Thus $P(\Delta t)$ is nearly proportional to
\[
    \frac12\cdot\frac1n\cdot2\Delta tg(X_1)\cdot(\Delta t/X_1)^{n-2}
\]
as $\Delta t \to 0$. It follows that the density of $P$ relative to ($n-1$)--volume, ie
\[
     \lim_{\Delta t\downarrow0}\frac{P(\Delta t)}{(\Delta t)^{n-1}V_{n-1}},
\]
is proportional to
\[
    \frac{g(X_1)}{X_1^{n-2}} = \frac{g(\|\vec X\|_\infty)}{\,\,\,\|\vec X\|_\infty^{n-2}}.
\]
We have, therefore, a G--R density with $f(s)=Kg(s)/s^{n-2}$, for some
constant $K$, and we must have (\ref{3.7}), ie $K=1/(n-1)V_{n-1}$,
because of the general relation between $f$ and $g$ for G--R
densities. 
\qed 
\end{proof}

To see whether a given density $g$ on $[0,1]$ generates a RBS distribution on $M(n)$ we must examine the values of $P\{X_k\in[a,b]\}$ in the model defined above. It is sufficient to compute $P\{X_1\in[0,t]\}$ since the coordinates $X_k$ have been symmetrized, ie they are interchangeable and $-X_k$ has the same distribution as $X_k$. Thus, the sampling vector $\vec X$ will have a RBS distribution on $M(n)$ exactly when 
\beq{3.10}{
       P\{X_1\in[0,t]\} = \frac{t}2\quad(t\in[0,1]). }
\end{equation}
Recalling the equidistribution of $Y_2,\dots,Y_n$, we compute:
\[ 
   P\{X_1\in[0,t]\} = P\{\pm=+\}P\{\sigma(1)=1\}P\{Y_1\in[0,t]\}+
\]\[
    P\{\pm=+\}P\{\sigma(1)\not=1\}P\{Y_2\in[0,t]\}+
         P\{\pm=-\}P\{\sigma(1)\not=1\}P\{Y_2\in[-t,0]\}
\]\[
    =\frac1{2n}\int_0^tg(s)\,ds + \frac12(1-\frac1n)P\{Y_1Z_2\in[-t,t]\}
\]\[
    =\frac1{2n}\int_0^tg(s)\,ds + \frac12(1-\frac1n)
                \int_0^1g(s)P\{Z_2\in[\frac{-t}s,\frac{t}s]\}\,ds.
\]
Let $P_n(r)=P\{|Z_2|\leq r\}$. Since $P_n(t/s)=1$ when $s<t$, our condition for a RBS distribution becomes: for all $t\in[0,1]$,
\beq{3.11}{
   \frac{t}2 = \frac12\int_0^tg(s)\,ds + \frac12(1-\frac1n)\int_t^1g(s)P_n(\frac{t}s)\,ds.
}\end{equation}
Differentiating with respect to $t$, we obtain the condition:
\beq{3.12}{
  1 \equiv \frac1ng(t)+(1-\frac1n)\int_t^1g(s)P_n'(\frac{t}s)\,\frac{ds}s\,\,\,(t\in[0,1]).
}\end{equation}

A little geometry reveals that $P_3(r)=r$ and $P_4(r)=r$, so that we may verify that (\ref{3.12}) is satisfied with $n=3,\,\,g_3(s)=3s^2$ and with $n=4,\,\,g_4(s)=4s^3$. This is one way to verify the discoveries of Robson and Gerow that $f_n(s)=C_ns$ yields RBS distributions for $n=3,4$. By a more involved geometric argument one may obtain $P_5(r)=(24s-s^3)/23$ and see that $g_5(s)=5s^4$ does not satisfy (\ref{3.12}) (with $n=5$), ie that the most ``natural'' generalization of the constructions of Robson and Gerow do not yield a RBS distribution for sample size 5. To use (\ref{3.12}) more systematically, we must first find some general expressions for $P_n'(r)$

Note that the ($n-2$)--volume on $\mathbb{R} \times [-1,1]^{n-2}
  \cap \{ \sum_2^n Z_k = -1 \}$ can be sampled uniformly by choosing
  $Z_3,\ldots Z_n$ independently and uniformly on $[-1,1]$ and then
  setting $Z_2 = -1 - \sum_3^n Z_k$.
  Thus the ($n-2$)--volume on $[-1,1]^{n-1}
  \cap \{ \sum_2^n Z_k = -1 \}$ can be sampled uniformly by repeating
  the above procedure until $Z_2 \in [-1,1]$. That is, until
  $\sum_3^nZ_k \in [0,2]$.
Thus, 
\beq{3.13}{
  P_n(s):= P\{|Z_2| \leq s\} = \frac{P\{\sum_3^nZ_k\in[1-s,1+s]\}}{P\{\sum_3^nZ_k\in[0,2]\}}\,\,\,(s\in[0,1]),
}\end{equation}
where the $Z_k$ are now independent and uniform over $[-1,1]$. Let $\phi$ be the uniform density over $[-1,1]$, ie $\phi=\frac12I_{[-1,1]}$. Then the density $\phi_n$ of $\sum_3^nZ_k$ is the ($n-2$)--fold convolution of $\phi$ with itself. Note that $\phi_n$ is an even function. It follows that, for $s\in[0,1]$,
\[
    P_n(s) = \frac{\int_{1-s}^{1+s}\phi_n(u)\,du}{\int_0^2\phi_n(u)\,du},
\]
and that
\[
   P_n'(s)=\frac{\phi_n(1+s)+\phi_n(1-s)}{\int_0^2\phi_n(u)\,du}
             =\frac{\phi_n(1+s)+\phi_n(1-s)}{2\phi_{n+1}(1)}.
\]
These functions, for $s\in[0,1]$, are polynomials in $s$. One way to see this, and to obtain explicit expressions for $P_n'(s)$, is to compute in terms of Laplace transforms. By the Laplace transform $\mathcal{L}\psi(s)$ of a function $\psi(t)$ (with left--bounded support) we mean
\[
   \mathcal{L}\psi(s) = \int_{-\infty}^\infty e^{-st}\psi(t)\,dt\quad(s>0),
\]
We shall use several well--known properties of $\mathcal{L}$; for example, $\mathcal{L}$ converts convolution products $\psi_1\star\psi_2$ into ordinary pointwise products:
\[
     \mathcal{L}\{\psi_1\star\psi_2\}(s)=\mathcal{L}\psi_1(s)\mathcal{L}\psi_2(s).
\]
Since $\mathcal{L}\phi(s)=(e^s-e^{-s})/2s$, it follows that 
\[
    \mathcal{L}\phi_n(s)=\frac{(e^s-e^{-s})^{n-2}}{2^{n-2}s^{n-2}}.
\]
Now $e^{rs}/s^{n-2}=\sL\{(t+r)_+^{n-3}/(n-3)!\}(s)$, where $t_+$ denotes $H(t)t$, $H(t)$ being the Heaviside function. Since $\sL$ is injective,
\beq{3.14}{
   \phi_n(t) = \frac1{(n-3)!2^{n-2}}\sum_{k=0}^{n-2}\binom{n-2}k (-1)^k (t+n-2-2k)_+^{n-3}.
}
\end{equation}
In evaluating $P_n'(s)$ (for $s\in[0,1]$), we need only apply (\ref{3.14}) for $t\in[0,2]$ and we obtain
\beq{3.15}{
   P_n'(s)=C_n\sum_{k=0}^{n-2}\binom{n-2}k(-1)^k\{(n-1-2k+s)_+^{n-3}+(n-1-2k-s)_+^{n-3}\},
}
\end{equation}
where $C_n=1/((n-3)!2^{n-1}\phi_{n+1}(1))$.
%  and $k^*=[(n-1)/2]$, since when $k>(n-1)/2$ we have $n-1-2k+s\leq-1+s\leq0$.

Let us consider the case of even $n$, say $n=2m$. 
Then in (\ref{3.15}), $2m-1-2k \pm s \leq 0$ if and only if $k\geq m$ 
(always with $s\in[0,1]$). So
\beq{3.16}{
   P_n'(s)=C_n\sum_{k=0}^{m-1}\binom{2(m-1)}k(-1)^k\{(2m-1-2k+s)^{n-3}+(2m-1-2k-s)^{n-3}\},
}
\end{equation}
%since the smallest value of $2m-1-2k\pm s$ is $2m-1-2(m-1)-1=0$. 
The odd case is the same but with one additional term.
In both cases, the odd powers of $s$ cancel, and the polynomial has the form
\beq{3.17}{
    P_n'(s) = C_n\sum_{j=0}^{\alpha_n}c_{j,n}s^{2j},
}
\end{equation}
for certain coefficients $c_{j,n}$, where $\alpha_n = \lfloor (n-3)/2 \rfloor$. 
These coefficients may be evaluated explicitly in any specific case, by reference to (\ref{3.15}); later we'll take a somewhat different point of view to derive some general properties of the $c_{j,n}$.

\begin{prop} \label{prop:3.3}
The G--R density corresponding to a probability density $g_n$ on [0,1] defines a RBS distribution on $M(n)$ exactly when
\beq{3.18}{
   \sL q_n(s)=\frac{n}{s(1+(n-1)\sL Q_n(s))},
}
\end{equation}
where $q_n(t)=H(t)g_n(e^{-t})$ and $Q_n(t)=H(t)P_n'(e^{-t})$.
\end{prop}

\begin{proof} 
The condition (\ref{3.12}), with $g=g_n$, may be rewritten in the form
\[
   n\equiv q_n(t)+(n-1)\int_{e^{-t}}^1 g_n(s)P_n'(\frac{e^{-t}}s)\,\frac{ds}s\quad(t\geq0).
\]
With the change of variable $s=e^{-u}$, this becomes
\[
    n\equiv q_n(t)+(n-1)\int_0^t q_n(u)P_n'(e^{-(t-u)})\,du\quad(t\geq0),
\]
ie
\[
   nH(t)=q_n(t)+(n-1)\int_{-\infty}^\infty q_n(u)Q_n(t-u)\,du,
\]
ie
\[
   nH = q_n + (n-1)q_n\star Q_n.
\]
Applying the Laplace transform we obtain
\[
   \frac{n}s = \sL q_n(s) + (n-1)\sL q_n(s)\sL Q_n(s).
\]
\qed
\end{proof}

Now in terms of (\ref{3.17}) we have
\[
    Q_n(t)=H(t)C_n\sum_{j=0}^{\alpha_n} c_{j,n}e^{-2jt},
\]
so that
\[
    \sL Q_n(s) = C_n\sum_{j=0}^{\alpha_n} \frac{c_{j,n}}{s+2j}.
\]
From Proposition~\ref{prop:3.3} it follows that, in order for $g_n$ to
generate a RBS distribution, we must have
\beq{3.19}{
   \sL q_n(s) = n/(s(1+(n-1)C_n\sum_{j=0}^{\alpha_n} \frac{c_{j,n}}{s+2j})).
}
\end{equation}

Using this condition one can check that the densities discovered by
Gerow and Robson are the only G--R densities for $n=3,4$.
Indeed evaluating~(\ref{3.19}), $\sL q_3(s) = \frac{3}{s+2}$ and $\sL q_4(s) = \frac{4}{s+3}$.
Thus, since $\sL$ is injective, $q_3(s) = 3e^{-2t}$ and $q_3(s) =
4e^{-3t}$.
So $g_3(t) = 3t^2$ and $g_4(t) = 4t^3$.

Continuing, using~(\ref{3.19}) one can check that $\sL q_5(s) =
  \frac{115(s+2)}{23s^2 + 130s + 192}$ and thus $q_5(t) =
  \frac{5}{191} e^{-\frac{65}{23}t}(191 \cos (\frac{\sqrt{191}}{23}t) -
  19 \sqrt{191} \sin ( \frac{\sqrt{191}}{23}t))$.
However $q_5(1.5)<0$. 
This contradicts the fact that $q(t) = g(e^{-t})$ and $g$ is a
  probability density on $[0,1]$.
Therefore there is no G--R density on $M(5)$.

We can generalize this non-existence result to higher $n$ without
explicitly calculating $\sL q_n$ and $q_n$.

Let 
\[
B_n(s) = \left(\prod_{j=0}^{\alpha_n} (s+2j) \right) \left( 1 +
  (n-1)C_n \sum_{j=0}^{\alpha_n} \frac{c_{j,n}}{s+2j} \right).
\]
Then $B_n(s)$ is a polynomial of degree $\alpha_n+1$ and 
\[
\sL q_n(s) = \frac{n \prod_{j=1}^{\alpha_n}(s+2j)}{B_n(s)}.
\]

Assume that $B_n(s)$ has $\alpha_n+1$ real distinct roots $a_0 > a_1 >
\ldots a_{\alpha_n}$.
Then for some $k\neq 0$ and some $\{b_j \neq 0\}$,
\[
\sL q_n(s) =
\frac{n\prod_{j=1}^{\alpha_n}(s+2j)}{k \prod_{j=0}^{\alpha_n}(s-a_j)} =
\sum_{j=0}^{\alpha_n} \frac{b_j}{s-a_j}.
\]
%for some $\{b_j\}$ where $b_j\neq 0$ for each $j$.
Since $\sL$ is injective,
\[
q_n(t) = \sum_{j=0}^{\alpha_n} b_j e^{a_j t}.
\]
Since $a_0 > a_1 > \ldots > a_{\alpha_n}$, for large values of $t$ the
sign of $q_n(t)$ equals the sign of $b_0$.
Since $g$ is a probability density, $q(t)\geq 0$ for all $t \geq 0$.
So $b_0 \geq 0$.
We have shown the following.

\begin{thm} \label{thm:3.4}
Let $n\geq 6$.
If $B_n(s)$ has $\alpha_n+1$ distinct real roots and $b_0<0$ then
there does not exist a G--R density on $M(n)$.
\end{thm}

With a few hours of computation on a desktop computer using a
computational program such as Maple one can verify that $B_n(s)$ has
$\alpha_n+1$ distinct real roots for $6 \leq n \leq 250$. 
We remark that we can show if $B_n(s)$ has $\alpha_n+1$
distinct real roots then $-3 < a_0 < -2$ and $b_0<0$.
However, to prove Theorem~\ref{thm:3.1} we do not need to use this fact.
The second condition can be verified computationally for $6 \leq n\leq
250$.
For example, for a given $n$ one can check that $-3 < a_0 < -2$.
Since for all $s \in (-3,-2)$, the sign of $B_n(s)$ is the opposite of
the sign of $\sL q_n(s)$, it follows that $b_0 < 0$ if and only if
$B_n(-3)<0$ and $B_n(-2)>0$. This proves Theorem~\ref{thm:3.1}.

%\bibliographystyle{elsart-num}
%\bibliography{my}

\begin{thebibliography}{10}
\expandafter\ifx\csname url\endcsname\relax
  \def\url#1{\texttt{#1}}\fi
\expandafter\ifx\csname urlprefix\endcsname\relax\def\urlprefix{URL }\fi

\bibitem{hammersleyMorton}
J.~M. Hammersley, K.~W. Morton, A new {M}onte {C}arlo technique: antithetic
  variates, Proc. Cambridge Philos. Soc. 52 (1956) 449--475.

\bibitem{arvidsenJohnsson}
N.~Arvidsen, T.~Johnsson, Variance reduction through negative correlation, a
  simulation study, J. Statist. Comput. Simulation (1982) 119--127.

\bibitem{hall:antithetic}
P.~Hall, Antithetic resampling for the bootstrap, Biometrika 76~(4) (1989)
  713--724.

\bibitem{craiuMeng:multiprocessPAC}
R.~V. Craiu, X.-L. Meng, Multiprocess parallel antithetic coupling for backward
  and forward {M}arkov {C}hain {M}onte {C}arlo, Annals of Statistics 33~(2) (2005)
  661--697.

\bibitem{craiuMeng:antitheticCouplingFPS}
R.~V. Craiu, X.-L. Meng, Antithetic coupling for perfect sampling, Bayesian
  Methods with Applications to Science, Policy, and Official Statistics -
  Selected Papers from ISBA 2000.

\bibitem{craiuMeng:chanceAndFractals}
R.~V. Craiu, X.-L. Meng, Chance and fractals, Chance 14 (2001) 47--52.

\bibitem{gerowMcCulloch}
K.~Gerow, C.~E. McCulloch, Simultaneously model-unbiased, design-unbiased
  estimation, Biometrics 56~(3) (2000) 873--878.

\bibitem{gerowHolbrook:mathIntelligencer}
K.~Gerow, J.~Holbrook, Statistical sampling and fractal distributions, Math.
  Intelligencer 18~(2) (1996) 12--22.

\bibitem{gerow:phdThesis}
K.~Gerow, Model-unbiased, unbiased-in-general estimation of a regression
  function, Ph.D. thesis, Cornell University (1993).

\bibitem{gerow:mastersThesis}
K.~Gerow, An algorithm for random balance sampling, Master's thesis, University
  of Guelph (1984).

\bibitem{royallCumberland}
R.~Royall, W.~Cumberland, An empirical study of the ratio estimator and
  estimators of its variance, Discussion paper in Journal of the American
  Statistical Society 76~(373) (1981) 66--88.

\bibitem{ramlochan:thesis}
R.~Ramlochan, Iterated function systems and fractal sampling, Master's thesis,
  University of Guelph (1990).

\bibitem{Hoel54}
P.~Hoel, Introduction to Mathematical Statistics, John Wiley \& Sons, 1954.

\bibitem{gerowHolbrook:constructionOfRBS}
K.~Gerow, J.~Holbrook, Construction of random balanced samples, working paper
  (1990).

\end{thebibliography}

\end{document}